\documentclass[journal]{IEEEtran}
\usepackage{color}
\usepackage{multirow}
\usepackage{cite}
\usepackage{balance}
\usepackage{amsmath}
\usepackage{enumitem}
\setlist[itemize]{noitemsep, topsep=0pt}

\usepackage{pdfsync}
\usepackage{amssymb}
\usepackage{bbold}
\usepackage{color}

\ifCLASSINFOpdf
\usepackage[pdftex]{graphicx}
\DeclareGraphicsExtensions{.pdf,.jpeg,.png, .jpg}
\else
\usepackage[dvips]{graphicx}

\DeclareGraphicsExtensions{.eps}
\fi
\graphicspath{{figs/}}
\usepackage{epstopdf}
\usepackage{booktabs}
\usepackage{amsmath}

\begin{document}

\title{Bi-Level Volt-VAR Optimization to Coordinate Smart Inverters with Voltage Control Devices}

\author{Rahul Ranjan Jha,~\IEEEmembership{Student Member,~IEEE,}
		Anamika Dubey,~\IEEEmembership{Member,~IEEE,}
        Chen-Ching Liu,~\IEEEmembership{Fellow,~IEEE, }
		and  Kevin P. Schneider,~\IEEEmembership{Senior Member,~IEEE}
      \thanks{R. Jha, A. Dubey, and C.C. Liu  are with the School of Electrical Engineering and Computer Science, Washington State University, Pullman, WA, 99164 e-mail: rahul.jha@wsu.edu,anamika.dubey@wsu.edu,liu@eecs.wsu.edu.}
        \thanks{K. P. Schneider is with the Pacific Northwest National Laboratory, Battelle Seattle Research Center, Seattle, WA 98109 USA e-mail: kevin.schneider@pnnl.gov.}
	}

\maketitle
\begin{abstract}
Conservation voltage reduction (CVR) uses Volt-VAR optimization (VVO) methods to reduce customer power demand by controlling feeder's voltage control devices. The objective of this paper is to present a VVO approach that controls system's legacy voltage control devices and coordinates their operation with smart inverter control. An optimal power flow (OPF) formulation is proposed by developing linear and nonlinear power flow approximations for a three-phase unbalanced electric power distribution system. A bi-level VVO approach is proposed where, Level-1 optimizes the control of legacy devices and smart inverters using a linear approximate three-phase power flow. In Level-2, the control parameters for smart inverters are adjusted to obtain an optimal and feasible solution by solving the approximate nonlinear OPF model. Level-1 is modeled as a Mixed Integer Linear Program (MILP) while Level-2 as a Nonlinear Program (NLP) with linear objective and quadratic constraints. The proposed approach is validated using 13-bus and 123-bus three-phase IEEE test feeders and a 329-bus three-phase PNNL taxonomy feeder. The results demonstrate the applicability of the framework in achieving the CVR objective. It is demonstrated that the proposed coordinated control approach help reduce feeder's power demand by reducing the bus voltages; the proposed approach maintains an average feeder voltage of 0.96 pu. A higher energy saving is reported during the minimum load conditions. The results and approximation steps are thoroughly validated using OpenDSS.
\end{abstract}
	

	%
	\IEEEpeerreviewmaketitle
		\vspace{-0.2 cm}
	\section*{Nomenclature}
	\addcontentsline{toc}{section}{Nomenclature}
	\begin{IEEEdescription}[\IEEEusemathlabelsep\IEEEsetlabelwidth{$V_1,V_2,V_3$}]
        \item [Sets]	
        \item[$\mathcal{G} = (\mathcal{N}, \mathcal{E})$] Directed graph for distribution system
        \item [$\mathcal{E}$] Set of distribution lines (branches) in $G$
        \item[$\mathcal{E_T}$] Set of branches with voltage regulator
        \item [$\mathcal{N}$] Set of buses (nodes) in $G$		
        \item[$\mathcal{N_C}$] Set of nodes with capacitor banks
		\item[$\mathcal{N_{DG}}$] Set of nodes with smart inverter connected DGs
        \item[$\Phi_i$] Set of phases of bus $i$ where, $\Phi_i \subseteq \{a,b,c\}$
        \item[$\Phi_{ij}$] $\{(pq): p\in \Phi_i, q\in \Phi_j\}$
        \item [Variables]
		\item[$(i,j)$] Line(branch) connecting nodes $i$ and $j$
        \item[$I_{ij}^p$] $I_{ij}^p = |I_{ij}^p|\angle{\delta_{ij}^p}$ is complex line current corresponding to phase $p\in\Phi_{i}$ where, $|I_{ij}^p|$ is magnitude and $\delta_{ij}^p$ is corresponding phase angle.
        \item[$p_{DG,i}^p$] Projected per-phase active power generated by $i^{th}$ DG at the current time instance.
		\item[$q_{DG,i}^p$] Available per-phase reactive power from $i^{th}$ DG at the current time instance.
        \item[$q_{cap,i}^{rated,p}$] Rated reactive power generated by capacitor bank connected to phase $p$ of bus $i$
        \item [$s$] Substation bus where, $s\in\mathcal{N}$
        \item[$s_{DG,i}^{rated,p}$] Rated per-phase apparent power capacity for DG connected to bus $i \in \mathcal{N_{DG}}$
        \item[$s_{L,i}^p$] $s_{L,i}^p = p_{L,i}^p + j q_{L,i}^p$ is complex power demand at bus $i$ corresponding to phase $p\in\Phi_i$ where, $p_{L,i}^p$ and $q_{L,i}^p$ are corresponding active and reactive power demand, respectively.
        \item[$S_{ij}^{pq}$] $S_{ij}^{pq} = P_{ij}^{pq} + j Q_{ij}^{pq}$ is complex power flow in branch $(i,j)$ corresponding to $(pq)\in\Phi_{ij}$, where, $P_{ij}^{pq}$ and $Q_{ij}^{pq}$ are corresponding active and reactive components, respectively.
        \item[$u_{tap,i}^p$] Binary control variable for voltage regulator tap position connected to phase $p$ of bus $i$
        \item[$u_{cap,i}^p$] Binary control variable for capacitor bank connected to phase $p$ of bus $i$  where, $i \in \mathcal{N_C}$
        \item[$V_i^p$] $V_{i}^p = |V_{i}^p|\angle{\theta_{i}^p}$ is complex voltage for $p\in\Phi_i$
        \item[$z_{ij}$] Complex three-phase impedance matrix for line $(i,j)\in \mathcal{E}$
		\item[$z_{ij}^{pq}$] $z_{ij}^{pq} = r_{ij}^{pq} + j x_{ij}^{pq}$ is an element of the complex impedance matrix $z_{ij}$ for branch $(i,j)$ where, $(pq)\in\Phi_{ij}$
	\end{IEEEdescription}

\vspace{-0.2 cm}
	
\section{Introduction}
\IEEEPARstart{C}{onservation} voltage reduction (CVR) is a technology to increase the energy efficiency of electric power distribution systems by reducing customer power demand through voltage control. The benefits of voltage control to energy savings are realized due to the sensitivity of customer loads to service voltages where decreasing the voltage helps reduce the demand \cite{schneider2010evaluation,EPRI}. Based on several studies and pilot projects, CVR can help achieve attractive energy savings. In fact, a recent study shows that CVR can help reduce the annual energy consumption by 3.04\% when implemented on all distribution feeders throughout the United States \cite{schneider2010evaluation}. 

Traditionally, CVR is accomplished by controlling feeder's legacy voltage control devices such as capacitor banks, load tap changers, and voltage regulators using Volt-VAR control (VVC) techniques. The feeder is operated at a lower service voltage range while still maintaining the service voltages within the recommended ANSI voltage limits (0.95 - 1.05 pu)\cite{ANSI}. In literature, several VVC methods have been proposed: 1) using autonomous or rule-based approach, 2) based on end-of-line measurements and 3) using integrated Volt-Var control (IVVC) based on real-time measurements \cite{baran1999volt,VVCBasicReview,VVCSanReview}. Several commercial VVC products are also available that perform IVVC function mostly using heuristic \cite{EPRIVVO}. Unfortunately, the available products only optimize the operation of legacy control devices. Recently, the integration of distributed generations (DGs) has increased in the distribution grid \cite{USRPS}. Most DGs are equipped with smart inverters that are capable of absorbing and supplying reactive power and thus controlling the feeder voltages locally that can help achieve additional CVR benefits \cite{Bokhari2}. Several researchers have worked on optimizing the reactive power dispatch from DGs and have proposed methods for smart inverter control using: 1) autonomous control, 2) distributed control, and 3) centralized control using optimal power flow (OPF) \cite{Farivar,Zhu1,Kekatos,Anese2,Su}.

The existing literature, however, presents several limitations. First, the available literature mostly fails to coordinate the control of feeder's legacy devices with smart inverters connected to DG. When mathematically modeling the VVO problem for both legacy and new devices, the optimization requires solving an OPF problem with both discrete and continuous variables. This results in a Mixed Integer Nonlinear Program (MINLP) that includes power flow equations for a three-phase unbalanced system. Some recent articles attempt to solve this problem \cite{VVO1,Dao,Wang1,Ranamuka,Hossan,Ding,liu2009reactive}. Unfortunately, these methods do not jointly optimize the control of legacy devices and new devices, do not apply to three-phase unbalanced system, or do not scale well even for mid-size system as the underlying problem is a MINLP. Second, the methods based on OPF assume a constant power load model thus fail to incorporate the voltage dependency of customer loads that is critical for modeling CVR effects. A voltage-dependent load model that can be easily incorporated within the optimization framework is called for. Third, the existing literature mostly solves a single-phase OPF problem. Distribution systems are largely unbalanced and require complete three-phase modeling to deliver a reasonable result. Recently, there has been some advancements in solving three-phase OPF problem, however, due to associated nonlinearities and mutual couplings, solving a three-phase OPF with only continuous decision variables is a challenging problem \cite{gan2014convex, Nanpeng}. Introducing discrete variables to three-phase OPF problem makes it even more challenging to solve. The aforementioned gap in literature calls for further research on enabling CVR for a modern distribution system.
The objective of this paper is to develop an OPF based bi-level approach for VVO to achieve CVR benefits for a three-phase unbalanced radial distribution system by simultaneously controlling both legacy devices and smart inverters. The proposed approach aims at addressing the aforementioned gaps in the literature and presents a scalable model for coordinated control of grid's all voltage control devices for CVR benefits. First, models for voltage-dependent loads and systems voltage control devices are developed to model CVR benefits, and valid formulations are proposed that can be easily incorporated within the three-phase OPF model. Next, we develop valid linear and nonlinear approximations for the three-phase power flow equations. Based on the proposed linear and nonlinear power flow approximations, a bi-level framework is proposed for CVR. Level-1 solves MILP to obtain set points for legacy device and smart inverters using the linear OPF model. Level-2 solves a NLP (with linear objective and quadratic constraints) based on nonlinear OPF to further optimize the smart inverter parameters and result in a feasible power flow solution.  Note that although energy savings are reported higher for meshed networks \cite{wang2018analysis}, a majority of distribution feeders in the United States are operated in radial topology. Therefore, this paper focuses on optimizing the operation of radially operated feeders. The specific contributions of this paper are detailed below.

\begin{itemize}[nolistsep,leftmargin=*]
   \item {\em Models for Voltage-dependent Loads and Voltage Control Devices:}  Mathematical models for voltage-dependent loads and grid's voltage control devices including capacitor banks, voltage regulators and smart inverters are proposed. For loads, a novel CVR-based load model is proposed that approximates the ZIP load model. The proposed models can be easily absorbed into both levels of optimization problems without changing the types of equations. 
  \item {\em Linear and Nonlinear Models for Three-phase Power Flow:} We develop valid linear and nonlinear approximations for three-phase unbalanced power flow model. The linear approximation is inspired by distflow equations but formulated for a three-phase unbalanced system. The nonlinear power flow is a new formulation and obtained by approximating the nonlinearities associated with the mutual coupling between the phases. Compared to the standard three-phase power flow formulations, the proposed model results in a reduced number of variables and introduces only nonlinearity of the nature of quadratic equality constraints.
  \item {\em Scalable Model for Coordinated Control of Legacy and New Voltage Control Devices for CVR:} The optimization problem for coordinated control of both legacy and new devices for a three-phase distribution feeder is a hard MINLP problem. To reduce complexity and enable scalability, we propose a bi-level approach by decomposing the MINLP into a MILP and a NLP. The scalability is demonstrated using IEEE 123-bus feeder (with 267 single-phase nodes) and 329-bus feeder (with 860 single-phase nodes). The proposed model solves 123-bus feeder within 4-mins and 329-bus feeder (after reduction) within 9-mins.
  \item {\em Validation using Multiple Test Feeders:} The proposed approach and all approximation steps are thoroughly validated against OpenDSS. First, the proposed three-phase approximate power flow  models are validated using IEEE 13-bus, IEEE 123-bus and 329-bus PNNL taxonomy feeders. Next, the accuracy of the proposed CVR-based load model is thoroughly validated against equivalent ZIP load models.
      Finally, the results for OPF are validated using OpenDSS. 
\end{itemize}

The rest of the paper is organized as follows. Section II presents the proposed approximate models three-phase power flow. Section III details mathematical models for distribution system equipment. Section IV presents the proposed bi-level VVO approach followed by results in Section V and conclusion in Section VI.

\vspace{-0.25cm}

\section{Three-Phase Unbalanced Electric Power Distribution System}

This section introduces the mathematical formulation for three-phase power flow based on branch-flow equations. Valid approximations are proposed to reduce the original formulation into a linear and an equivalent quadratic formulation.

\vspace{-0.5cm}
\subsection{Three-Phase Power Flow using Branch Flow Model}
Let, there be directed graph $\mathcal{G} = (\mathcal{N}, \mathcal{E})$ where $\mathcal{N}$ denotes set of buses and $\mathcal{E}$ denotes set of lines. Each line connects ordered pair of buses $(i,j)$ between two adjacent nodes $i$ and $j$. Let, $\{a,b,c\}$ denotes the three phases of the system and $\Phi_i$ denotes set of phases on bus $i$. For each bus $i \in \mathcal{N}$, let, phase $p$ complex voltage is given by $V_i^{p}$ and phase $p$ complex power demand is $s_{L,i}^p$. Let, $V_i := [V_i^{p}]_{p \in \Phi_i}$ and $s_{L,i} := [s_{L,i}^{p}]_{p \in \Phi_i}$. For each line, let, $p$ phase current be $I_{ij}^{p}$ and define, $I_{ij} := [I_{ij}^{p}]_{p \in \Phi_i}$. Let, $z_{ij}$ be the phase impedance matrix.

\vspace{-0.5cm}
\begin{eqnarray}
v_j  = v_i - (S_{ij}z_{ij}^H+z_{ij}S_{ij}^H) + z_{ij}l_{ij} z_{ij}^H \\
\text{diag}(S_{ij} - z_{ij}l_{ij}) =  \sum_{k:j \rightarrow k}{\text{diag}(S_{jk})}  + s_{L,j}  \\
\left[
  \begin{array}{cc}
    v_i & S_{ij} \\
    S_{ij}^H & l_{ij} \\
  \end{array}
\right] = \left[
  \begin{array}{c}
    V_i\\
    I_{ij}\\
  \end{array}
\right]
\left[
  \begin{array}{cc}
    V_i\\
    I_{ij}\\
  \end{array}
\right]^H \\
\left[
  \begin{array}{cc}
    v_i & S_{ij} \\
    S_{ij}^H & l_{ij} \\
  \end{array}
\right] : - \text{Rank-1 PSD Matrix}
\end{eqnarray}

A three-phase power flow formulation for a radial system based on branch flow relationship is given in \cite{gan2014convex} and detailed in (1)-(4). Here, (1) represents voltage drop equation, (2) corresponds to power balance equation, (3) are variable definitions for power flow quantities, and (4) is a Rank-1 constraint that makes the associated optimization problem non-convex. In the literature, methods are proposed to obtain a relaxed convex problem \cite{Low1,Low2}, however, it is difficult to obtain a feasible solution from relaxed problem for a three-phase system \cite{Nanpeng}. Moreover, it is difficult to extend the power flow model detailed in (1)-(4) for voltage-dependent loads and system's legacy control devices. A new three-phase power flow model for OPF problem is called for that can easily incorporate system's critical components while not significantly increasing the inherent nonlinearity.

\vspace{-0.5cm}

\subsection{Approximate Three-Phase Power Flow Equations}
In this section, we present valid linear and nonlinear power flow equations by approximating (1)-(4). Fundamentally, there are two reasons for nonlinearity in power flow equations: nonlinear relationship between power, voltage, and currents, and mutual coupling in a three-phase system. In the proposed formulation, the nonlinearity resulting from mutual coupling between the three phases is approximated. A phase-decoupled formulation by decoupling the branch power flow and voltage equations on a per-phase basis is obtained. The resulting three-phase power flow model characterizes the power flow equations using a fewer number of variables.

Define: $v_i^p = (V_i^p)^2 \text{ where, } p \in  \Phi_{i}$, $l_{ij}^{pq} = (|I_{ij}^p|\times |I_{ij}^q|) \text{ where, } (pq) \in  \Phi_{ij}$, $\delta_{ij}^{pq} = \delta_{ij}^{p} - \delta_{ij}^{q}$, $S_{ij}^{pq} = P_{ij}^{pq} + j Q_{ij}^{pq}, \text{ where, } (pq) \in  \Phi_{ij}$, and $z_{i,j}^{pq} = r_{ij}^{pq} + j x_{ij}^{pq}, \text{ where, } (pq) \in  \Phi_{ij}$. Note that, $ (l_{ij}^{pq})^2 = l_{ij}^{pp} \times l_{ij}^{qq}$ and $\delta_{ij}^{pq}$ is angle difference between branch currents $I_{ij}^p$ and $I_{ij}^q$.

\subsubsection{Assumption 1 - Approximating Phase Voltages} The phase voltages are assumed to exactly $120^\circ$ degree apart. Moreover, it is assumed that the degree of unbalance in voltage magnitude is not large. Both conditions are valid for a distribution system as specified by the ANSI limits for bus voltages and phase unbalance \cite{ANSI}. 

\vspace{-0.2cm}
\begin{equation}\label{}
  \dfrac{V_i^a}{V_i^b} \simeq \dfrac{V_i^b}{V_i^c} \simeq \dfrac{V_i^c}{V_i^a} = e^{j*2\pi/3} = a^2
\end{equation}
Further note that
\begin{equation}\label{}
  S_{ij}^{pq} = V_{i}^{p} \times I_{ij}^{q} \hspace{6 pt} \text{and} \hspace{6 pt} S_{ij}^{qq} = V_{i}^{q} \times I_{ij}^{q}
\end{equation}
Using (5) and (6), we express $S_{ij}^{pq}$ as a function of $S_{ij}^{qq}$ i.e. $S_{ij}^{pq} = \dfrac{V_i^p}{V_i^q}\times S_{ij}^{qq}$, where, $\dfrac{V_i^p}{V_i^q}$ is a constant (5). The off-diagonal elements of $S_{ij}$ are approximated using diagonal terms which help reduce the number of power flow variables.

Note that the above conditions do not imply that a single-phase power flow will be sufficient to represent a distribution system. First, an equivalent single-phase model cannot represent two-phase or single-phase lines and loads. Second, it is imperative to solve for an unbalanced power flow even though the degree of voltage unbalance is less. 

\subsubsection{Assumption 2 - Approximating Angle Difference between Phase Currents}
On expanding (1) and (2), nonlinearities are introduced as trigonometric functions of angle difference between the phase currents on a given three or two-phase line. Let, for a given line $(i,j)$, the phase currents for phases $p$ and $q$ be $I_{ij}^p = |I_{ij}^p| \angle{\delta_{ij}^p}$, and $I_{ij}^q = |I_{ij}^q| \angle{\delta_{ij}^q}$. Then the angle difference between the phase currents of a given bus $i$ i.e. $\delta_{ij}^{pq} = \angle{\delta_{ij}^p} - \angle{\delta_{ij}^q}$. We observe terms corresponding to $\sin(\delta_{ij}^{pq})$ and $\cos(\delta_{ij}^{pq})$ in power flow expressions. These terms significantly increase the complexity of the OPF problem.

In the proposed formulation, the phase angle difference between branch currents are approximated and modeled as a constant variable. An approximate value of $\delta_{ij}^{pq}$ for branch $(i,j)$ is calculated by solving an equivalent distribution power flow with system loads modeled as constant impedance loads. We assume $\delta_{ij}^{pq}$ to be constant and equal to the one obtained by solving power flow with constant impedance load model.  Note that the constant impedance load model is only used to approximate $\delta_{ij}^{pq}$. This assumption does not limit the type of load that can be incorporated in the proposed OPF model; it can easily incorporate all load types including constant power and constant current loads, as detailed in Section IV.

\subsubsection{Power Flow Equations}
The power flow equations defined in (1)-(4) are expanded. Using the above approximations, we are able to redefine the power flow equations in (1)-(4) as a set of linear and non-linear equations shown in (7)-(11). Here, (7)-(9) are linear in $v_i^p$, $l_{ij}^{pq}$, and $S_{ij}^{pq}$. Note that the total number of variables in the proposed formulation are $15\times(n-1)$, where $n$ is the number of nodes; the original formulation (1)-(4) had a total of $36 \times(n-1)$ variables.
\vspace{-0.2cm}
\begin{eqnarray}
  \nonumber &P_{ij}^{pp} - \sum_{q \in \Phi_j}{l_{ij}^{pq}\left(r_{ij}^{pq} \cos(\delta_{ij}^{pq})- x_{ij}^{pq} \sin(\delta_{ij}^{pq})\right)} \\ &= \sum_{k:j \rightarrow k}P_{jk}^{pp} + p_{L,j}^p\\
  \nonumber &Q_{ij}^{pp} - \sum_{q \in \Phi_j}{l_{ij}^{pq}\left(x_{ij}^{pq} \cos(\delta_{ij}^{pq})+ r_{ij}^{pq} \sin(\delta_{ij}^{pq})\right)}\\&= \sum_{k:j \rightarrow k}Q_{jk}^{pp} + q_{L,j}^p
\end{eqnarray}
\begin{eqnarray}
 \nonumber &v_j^p = v_i^p - \sum_{q \in \Phi_j}{2 \mathbb{Re}\left[S_{ij}^{pq} (z_{ij}^{pq})^*\right]} + \sum_{q \in \Phi_j}{z_{ij}^{pq} l_{ij}^{qq}} \\
  &+ \sum_{q1,q2 \in \Phi_j, q1 \neq q2}{2\mathbb{Re}\left[ z_{ij}^{pq1} l_{ij}^{q1q2}\left(\angle(\delta_{ij}^{q1q2})\right)(z_{ij}^{pq2})^*\right]}
\end{eqnarray}
\begin{eqnarray}
  (P_{ij}^{pp})^2 + (Q_{ij}^{pp})^2 = v_i^p  l_{ij}^{pp}\\
  (l_{ij}^{pq})^2 = l_{ij}^{pp}  l_{ij}^{qq}
\end{eqnarray}

\begin{itemize}[leftmargin=*]
  \item (7) is written for all $(i,j)\in \mathcal{E}$ and represents the equation for active power flow on branch $(i,j)$ for phase $pp \in \Phi_{ij}$ . 
      Since $\cos(\delta_{ij}^{pq})$and $\sin(\delta_{ij}^{pq})$ are assumed to be constant, this equation is linear in $P_{ij}^{pp}$, $P_{jk}^{pp}$, and $l_{ij}^{pq}$.
  \item (8) is written for all $(i,j)\in \mathcal{E}$ and represents the equation for reactive power flow on branch $(i,j)$ for phase $pp \in \Phi_{ij}$. 
      Similar to (7), (8) is linear in $Q_{ij}^{pp}$, $Q_{jk}^{pp}$, $l_{ij}^{pq}$.
  \item (9) represents the equation for voltage drop between the two nodes of branch $(i,j)$ corresponding to phase $p$. The equation is also linear in $v_i^p$, $v_j^p$, and $l_{ij}^{pq}$.
  \item (10) relates per phase complex power flow in branch $(i,j)$ to phase voltage and phase currents. This is a non-linear quadratic equality constraint.
  \item (11) simply relates current variables previously defined, i.e. $l_{ij}^{pq} = (|I_{ij}^p|\times |I_{ij}^q|), pq \in  \Phi_{ij}$. 
\end{itemize}

\vspace{-0.4cm}
\subsection{Linear Three-Phase AC Power Flow Approximation}
The linear approximation assumes the branch power loss are relatively smaller as compared to the branch power flow \cite{gan2014convex}. The impact of power loss on active and reactive power branch flow equations and on voltage drop equations is ignored. After approximating (7)-(11), we obtain linearized AC branch flow equations as shown in (12)-(13). Here (12) corresponds to linearized active and reactive power flow and (13) corresponds to voltage drop equations. 
\vspace{-0.2cm}
\begin{eqnarray}
  P_{ij}^{pp} = \sum_{k:j \rightarrow k}P_{jk}^{pp} + p_{L,j}^p \hspace{0.1cm} \text{and}  \hspace{0.1cm}
  Q_{ij}^{pp} = \sum_{k:j \rightarrow k}Q_{jk}^{pp} + q_{L,j}^p \\
 v_j^p = v_i^p - \sum_{q \in \Phi_j}{2 \mathbb{Re}\left[S_{ij}^{pq} (z_{ij}^{pq})^*\right]} \forall j \in Y_i
\end{eqnarray}

The AC linearized power flow is significantly accurate in representing bus voltages. 
The linearized AC power flow, although does not include the impact of power loss on voltage drop, it does incorporate the impact of power flow due to load. Since power losses are significantly small as compared to the branch flow due to load demand, the obtained feeder voltages are good approximation of the actual feeder voltages \cite{gan2014convex}.

\vspace{-0.1cm}

\section{Distribution System Equipment Models}

This section details the models for capacitor banks, voltage regulators, smart inverters and voltage-dependent customer loads. 
The approximate power flow equations developed in Section II are a function of $v_i^p = |V_i^p|^2$. The equipment models are, therefore, parameterized based on respective control variables and the $v_i^p$. A new CVR based load model is developed to represent the power demand as a function of $v_i^p$. The ZIP coefficients for the load are used to obtain equivalent CVR coefficients. 
Note that the equipment and load models proposed in this section are specifically designed so that they can be easily absorbed within the approximate power flow equations defined in (7)-(13) without changing their type.

\vspace{-0.5cm}
\subsection{Voltage Regulator}
A 32-step voltage regulator with a voltage regulation range of $\pm10\%$ is assumed. The series and shunt impedance of the voltage regulator are ignored as these have very small value \cite{kers}. Let, $a^p$ be the turn ratio for the voltage regulator connected to phase $p$ of line $(i,j)$. Then $a^p$ can take values between 0.9 to 1.1 with each step resulting in a change of 0.00625 pu. An additional node $i'$ is introduced to model the current equations. The control for regulator is defined using binary variables. Let, for $u_{tap,i}^p \in \{0,1\}$ be a binary variable defined for each regulator step position i.e. $i \in (1,2,...,32)$. Also define a vector $b_i \in \{0.9, 0.90625 , ..., 1.1\}$. Then $V_i^p$, $V_j^p$, $I_{ii'}^p$, and $I_{i'j}^p$ where $p \in \Phi_i \cap \Phi_j$ are given as follows:
\begin{eqnarray} \label{eq1}
V_j^p = V_{i'}^p = a^p  V_i^p \hspace{0.1cm} \text{and} \hspace{0.1cm}
I_{ii'}^p = a^p  I_{i'j}^p
\end{eqnarray}
where, $a^p = \sum\limits_{i=1}^{32} b_i u_{tap,i}^p$ and  $\sum\limits_{i=1}^{32} u_{tap,i}^p = 1$.

In order to express (14) as a function of $v_i^p = (V_i^p)^2$, $v_j^p = (V_j^p)^2$, $l_{ii'}^{pp}=(I_{ii'}^p)^2$, and $l_{i'j}^{pp}=(I_{i'j}^p)^2$ we take square of (14) and define $a_p^2 = A_p$ and $b_i^2 = B_i$. Further realizing that $(u_{tap,i}^p)^2 = u_{tap,i}^p$, (14) can be reformulated as (15). 
\begin{eqnarray} \label{eq1reg}
v_j^p = A^p \times v_i^p \hspace{0.2cm} \text{and} \hspace{0.2cm}
l_{ii'}^{pp} = A^p  l_{i'j}^{pp}
\end{eqnarray}

\vspace{-0.7cm}
\subsection{Capacitor Banks}
The per-phase model for capacitor banks is developed. The reactive power generated by capacitor bank, $q_{cap,i}^{p}$, is defined as a function of binary control variable $u_{cap,i}^p \in \{0,1\}$ indicating the status (On/Off) of the capacitor bank, its rated per-phase reactive power $q_{cap,i}^{rated,p}$, and the square of the bus voltage at bus $i$ for phase $p$, $v_{i}^p$.
\begin{equation} \label{eq2}
q_{cap,i}^{p} = u_{cap,i}^p q_{cap,i}^{rated,p} v_{i}^p
\end{equation}

The capacitor bank model is assumed to be voltage dependent and provides reactive power as a function of $v_{i}^p$ when connected, i.e. $u_{cap,i}=1$. For a three-phase capacitor bank, a common control variable, $u_{cap,i}^p$, is defined for each phase.

\vspace{-0.5cm}
\subsection{Distributed Generation with Smart Inverters}
A per-phase model for reactive power support from smart inverter connected to DGs is developed. The DGs are modeled as negative loads with a known active power generation equal to the forecasted value. The reactive power support from DG depend upon the rating of the smart inverter. Let, the rated per-phase apparent power capacity for smart inverter connected to $i^{th}$ DG be $s_{DG,i}^{rated,p}$ and the forecasted active power generation be $p_{DG,i}^p$. The available reactive power, $q_{DG,i}^p$ from the smart inverter is given by (17) which is a box constraint.
\begin{equation} \label{eq3}
\small
 -\sqrt{(s_{DG,i}^{rated,p})^2 - (p_{DG,i}^p)^2} \leq q_{DG,i}^p \leq \sqrt{(s_{DG,i}^{rated,p})^2 - (p_{DG,i}^p)^2}
\end{equation}

\vspace{-0.5cm}
\subsection{Voltage-Dependent Model for Customer Loads}

The most widely acceptable load model is the ZIP model which is a combination of constant impedance (Z), constant current (I) and constant power (P)) characteristics of the load \cite{Bokhari}. The mathematical representation of the ZIP model for the load connected at phase $p$ of bus $i$ is given by (18)-(19).
\begin{eqnarray} \label{eq4}
p_{L,i}^p = p_{i,0}^p \left[k_{p,1} \left(\dfrac{V_i^p}{V_0}\right)^2 + k_{p,2} \left(\dfrac{V_i^p}{V_0}\right)+ k_{p,3}\right]\\
q_{L,i}^p = q_{i,0}^p \left[k_{q,1} \left(\dfrac{V_i^p}{V_0}\right)^2 + k_{q,2} \left(\dfrac{V_i^p}{V_0}\right)+ k_{q,3}\right]
\end{eqnarray}
where, $k_{p,1}+k_{p,2}+k_{p,3} = 1$, $k_{q,1}+k_{q,2}+k_{q,3} = 1$, $p_{i,0}^p$ and $q_{i,0}^p$ are per-phase load consumption at nominal voltage, $V_0$.

The ZIP load model represented in (18)-(19) are a function of both $V_i^p$ and $v_i^p = (V_i^p)^2$. Including (18) and (19) to OPF formulation will make (7),(8),(13), and (14), earlier linear in $v_i^p$, nonlinear. Here we develop an equivalent load model for voltage-dependent loads using the definition of CVR factor. Next, an equivalence between ZIP parameters and proposed CVR factors is obtained.

CVR factor is defined as the ratio of percentage reduction in active or reactive power to the percentage reduction in bus voltage. Let CVR factor for active and reactive power reduction be $CVR_{p}$, and $CVR_{q}$, respectively defined in (25).
\begin{eqnarray} \label{eq7}
CVR_{p} = \dfrac{d p_{L,i}^p}{p_{i,0}^p} \dfrac{V_0}{dV_i^p} \hspace{0.2cm} \text{and} \hspace{0.2cm} CVR_{q} = \dfrac{d q_{L,i}^p}{q_{i,0}^p} \dfrac{V_0}{dV_i^p}
\end{eqnarray}
where, $p_{L,i}^p = p_{i,0}^p + d p_i^p$ and $q_{L,i}^p = q_{i,0}^p + d q_i^p$.
Furthermore, $v_i^p = (V_i^p)^2$. Therefore, $d v_i^p =  2V_i^p dV_i^p$. Assuming $V_i^p \approx V_0$ and $d v_i^p = v_i^p-(V_0)^2$, we obtain:
\begin{eqnarray}
p_{L,i}^p = p_{i,0}^p + CVR_{p}\dfrac{p_{i,0}^p}{2}\left(\dfrac{v_i^p}{V_0^2}-1\right)\\
q_{L,i}^p = q_{i,0}^p + CVR_{q}\dfrac{q_{i,0}^p}{2}\left(\dfrac{v_i^p}{V_0^2}-1\right)
\end{eqnarray}
Note that the CVR based load model detailed in (21) and (22) is linear in $v_i^p$, thus can be easily included in approximate power flow equations (7)-(13).
The CVR factors, $CVR_{p}$ and $CVR_{q}$ are estimated from the ZIP coefficients of the load. On differentiating the ZIP model detailed in (18) and (19) and assuming $V_0 = 1$ p.u., we obtain:
\begin{eqnarray} \label{eq10}
\dfrac{d p_{L,i}^p}{dV_i^p} = p_{i,0}^p \left(2 k_{p,1} V_i^p +  k_{p,2}\right)\\
\dfrac{d q_{L,i}^p}{dV_i^p} = q_{i,0}^p \left(2 k_{q,1} V_i^p +  k_{q,2}\right)
\end{eqnarray}
Using (20), (23), (24) and assuming $V_i^p \approx V_0$, we obtain (25). Using (25), the CVR factors for customer loads can be obtained from the ZIP coefficients.
\begin{equation} \label{eq12}
CVR_p = 2 k_{p,1} + k_{p,2} \hspace{0.1cm} \text{and} \hspace{0.1cm} CVR_q = 2 k_{q,1} + k_{q,2}
\end{equation}
\section{Proposed Bi-level Volt-VAR Optimization}
The primary function of VVO is to use voltage control to 
1) reduce energy consumption, 2) reduce system losses, and 3) regulate feeder voltages. The problem of coordinating the control of system's legacy devices and smart inverters results in a MINLP problem. To reduce complexity and ensure scalability, a bi-level approach is proposed.
\begin{enumerate}
\item	Level 1: Develops a 15-min schedule for legacy devices and smart inverter reactive power demand set-points with the objective of minimizing the active power consumption for the feeder based on a MILP formulation.

\item	Level 2: Develops revised 15-min schedule for smart inverter controls using a NLP formulation to achieve feasible three-phase power flow solutions. Level-2 uses the nonlinear power flow formulation proposed in Section II.C and obtains revised set points for smart inverter control that ensure feasible power flow solutions. 
\end{enumerate}

\vspace{-0.3cm}

\subsection{Level 1 - MILP Formulation for Coordinated Control}
The objective of this stage is to minimize the power consumption for the feeder by controlling voltage regulators, capacitor banks, and smart inverters while ensuring that the voltage limits are satisfied. The control of legacy devices introduces integer variables into the optimization problem. A linear three-phase AC power flow is used and resulting problem is a MILP formulation as detailed in (26)-(37). The objective is to minimize the sum of three-phase active power flowing out of the substation bus at time $t$ (26). Here, $s \in \mathcal{N}$ denotes the substation bus. Since, the distribution feeder is radial, the substation power equals net feeder power demand.

\vspace{0.2cm}
\begin{minipage}{23.5em}
		\flushleft
		\small
		Variables: \\ $u_{tap,i}^p(t)$, $u_{cap,i}^p(t)$, $q_{DG,i}^p(t)$, $v_i^p(t)$, $P_{ij}^{pp}(t)$, $Q_{ij}^{pp}(t)$, $S_{ij}^{pq}(t)$
		\begin{flalign}\label{eq15}
		\text{Minimize:}  \ \ \sum_{p\in \Phi_s,j:s \rightarrow j}{P_{sj}^p(t)}		&&
		\end{flalign}
		Subject to:
        \begin{flalign}\label{eq17}
            &P_{ij}^{pp}(t) = \sum_{k:j \rightarrow k}P_{jk}^{pp}(t) + p_{L,j}^p(t)- p_{DG,i}^p(t) \hspace{0.1cm} \forall i\in\mathcal{N} \\
            &Q_{ij}^{pp}(t) = \sum_{k:j \rightarrow k}Q_{jk}^{pp}(t) + q_{L,j}^p(t)- q_{DG,i}^p(t) - q_{C,i}^p \hspace{0.1cm} \forall i\in\mathcal{N} \\
            &v_j^p(t) = v_i^p(t) - \sum_{q \in \Phi_j}{2 \mathbb{Re}\left[S_{ij}^{pq}(t)(z_{ij}^{pq})^*\right]} \hspace{0.1cm} \forall j \in Y_i \\
            &p_{L,i}^p(t) = p_{i,0}^p(t) + CVR_{p}(t)\dfrac{p_{i,0}^p(t)}{2}(v_i^p(t)-1)  \forall i\in\mathcal{N_L} \\
            &q_{L,i}^p(t) = q_{i,0}^p(t) + CVR_{q}(t)\dfrac{q_{i,0}^p(t)}{2}(v_i^p(t)-1)  \forall i\in\mathcal{N_L} \\
            &v_j^p(t) = A_i^p(t) v_i^p(t) \forall (i,j) \in \mathcal{E_T} \\
            &A_i^p(t) = \sum\limits_{i=1}^{32} B_i u_{tap,i}^p(t), \sum\limits_{i=1}^{32} u_{tap,i}^p(t) = 1 \forall (i,j) \in \mathcal{E_T}  \\
            &q_{C,i}^{p}(t) = u_{cap,i}^p(t)  q_{cap,i}^{rated,p} v_i^p(t)  \hspace{0.1cm} \forall (i) \in \mathcal{N_C} \\
            &q_{DG,i}^p(t) \leq \sqrt{(s_{DG,i}^{rated,p})^2 - (p_{DG,i}^p)^2(t)} \hspace{0.2cm} \forall (i) \in \mathcal{N_{DG}} \\
            & q_{DG,i}^p(t) \geq -\sqrt{(s_{DG,i}^{rated,p})^2 - (p_{DG,i}^p)(t)^2} \hspace{0.2cm} \forall (i) \in \mathcal{N_{DG}}\\
            &(V_{min})^2\leq v_i^p(t) \leq (V_{max})^2 \hspace{12pt} \forall i\in \mathcal{N}
        \end{flalign}
		\end{minipage}
\vspace{0.2cm}
\begin{itemize}[noitemsep,topsep=0pt,leftmargin=*]
  \item Constraints (27)-(29) are linear AC power flow constraints.
  \item Constraints (30)-(31) define CVR based load model.
  \item Constraints (32)-(33) define regulator control equations.
  \item Constraint (34) defines equations for capacitor control.
  \item Constraints (35)-(36) define control equations for reactive power dispatch at time $t$ from smart inverters.
  \item Constraints (37) defines operating limits for feeder voltages. 
\end{itemize}

\vspace{-0.2cm}

\subsection{Level 2 - NLP Problem for Smart Inverter Control}
Level-1 uses a linear three-phase power flow model that approximates the losses. The solutions although feasible for linear power flow formulation, may violate the critical operating constraints of the feeder. The objective of this stage is to adjust the set-points of smart inverter control variables in order to obtain an optimal and feasible three-phase nonlinear power flow solution. The discrete control variables, $u_{tap,i}^p(t)$, $u_{cap,i}^p(t)$, are assumed to be fixed as obtained in Level-1. The optimal control set points for reactive power dispatch from smart inverters are obtained by solving the NLP problem (with linear objective and quadratic constraints) defined in (38)-(50).

\vspace{0.2cm}
\begin{minipage}{23.5em}
		\flushleft
		\small
		Variables: $q_{DG,i}^p(t)$, $v_i^p(t)$, $P_{ij}^{pp}(t)$, $Q_{ij}^{pp}(t)$, $S_{ij}^{pq}(t)$, $l_{ij}^{pq}(t)$
		\begin{flalign}\label{eq15nl}
		\text{Minimize:}  \ \ \sum_{p\in \Phi_s,j:s \rightarrow j}{P_{sj}^p(t)}		&&
		\end{flalign}
		Subject to:
        \begin{flalign}\label{eq17nl}
        \nonumber &P_{ij}^{pp}(t) - \sum_{q \in \Phi_j}{l_{ij}^{pq}(t)\left(r_{ij}^{pq} \cos(\delta_{ij}^{pq}(t))- x_{ij}^{pq} \sin(\delta_{ij}^{pq}(t))\right)} \\ &= \sum_{k:j \rightarrow k}P_{jk}^{pp}(t) + p_{L,j}^p(t)- p_{DG,i}^p(t) \hspace{0.1cm} \forall i\in\mathcal{N}\\
        \nonumber &Q_{ij}^{pp}(t) - \sum_{q \in \Phi_j}{l_{ij}^{pq}(t)\left(x_{ij}^{pq} \cos(\delta_{ij}^{pq}(t))+ r_{ij}^{pq} \sin(\delta_{ij}^{pq}(t))\right)}\\&= \sum_{k:j \rightarrow k}Q_{jk}^{pp}(t) - q_{DG,i}^p(t) - q_{C,i}^p(t) \hspace{0.1cm} \forall i\in\mathcal{N} \\
        \nonumber &v_j^p(t) = v_i^p(t) - \sum_{q \in \Phi_j}{2 \mathbb{Re}\left[S_{ij}^{pq}(t) (z_{ij}^{pq})^*\right]} + \sum_{q \in \Phi_j}{z_{ij}^{pq} l_{ij}^{qq}(t)} \\
        &+ \sum_{q1,q2 \in \Phi_j, q1 \neq q2}{2\mathbb{Re}\left[ z_{ij}^{pq1} l_{ij}^{q1q2}(t)\left(\angle(\delta_{ij}^{q1q2}(t))\right)(z_{ij}^{pq2})^*\right]} \\
        &(P_{ij}^{pp}(t))^2 + (Q_{ij}^{pp}(t))^2 = v_i^p(t)  l_{ij}^{pp}(t) \hspace{0.1cm} \forall (i,j)\in\mathcal{E}\\
        &(l_{ij}^{pq}(t))^2 = l_{ij}^{pp}(t)  l_{ij}^{qq}(t) \hspace{0.1cm} \forall (i,j)\in\mathcal{E}\\
        &p_{L,i}^p(t) = p_{i,0}^p(t) + CVR_{p}(t)\dfrac{p_{i,0}^p(t)}{2}(v_i^p(t)-1)  \forall i\in\mathcal{N_L} \\
        &q_{L,i}^p(t) = q_{i,0}^p(t) + CVR_{q}(t)\dfrac{q_{i,0}^p(t)}{2}(v_i^p(t)-1)  \forall i\in\mathcal{N_L} \\
        &v_j^p(t) = A_i^p(t) v_i^p(t), A_i^p(t) = \sum\limits_{i=1}^{32} B_i u_{tap,i}^p(t) \forall (i,j) \in \mathcal{E_T} \\
        &q_{C,i}^{p}(t) = u_{cap,i}^p(t) q_{cap,i}^{rated,p} v_i^p(t)  \hspace{0.1cm} \forall (i) \in \mathcal{N_C} \\
        &q_{DG,i}^p(t) \leq \sqrt{(s_{DG,i}^{rated,p})^2 - (p_{DG,i}^p)^2(t)} \hspace{0.2cm} \forall (i) \in \mathcal{N_{DG}} \\
        &q_{DG,i}^p(t) \geq -\sqrt{(s_{DG,i}^{rated,p})^2 - (p_{DG,i}^p)(t)^2} \hspace{0.2cm} \forall (i) \in \mathcal{N_{DG}}\\
        &(V_{min})^2\leq v_i^p(t) \leq (V_{max})^2 \hspace{12pt} \forall i\in \mathcal{N}
        \end{flalign}
		\end{minipage}

\vspace{0.2cm}
\begin{itemize}[noitemsep,topsep=0pt,leftmargin=*]
  \item Constraints (39)-(43) are approximate nonlinear AC power flow equations defined at time $t$.
  \item Constraints (44)-(45) define CVR based load model.
  \item Constraints (46) define regulator control equations. Note that $u_{tap,i}^p$ is known from Level-1 solution.
  \item Constraint (47) defines control equations for capacitor banks at time $t$. Note that $u_{cap,i}^p$ is known from Level-1 solution.
  \item Constraints (48)-(49) define control equations for reactive power dispatch at time $t$ from smart inverters.
  \item Constraints (50) defines operating limits for feeder voltages. 
\end{itemize}

\begin{figure}[t]
\includegraphics[width=3.4in]{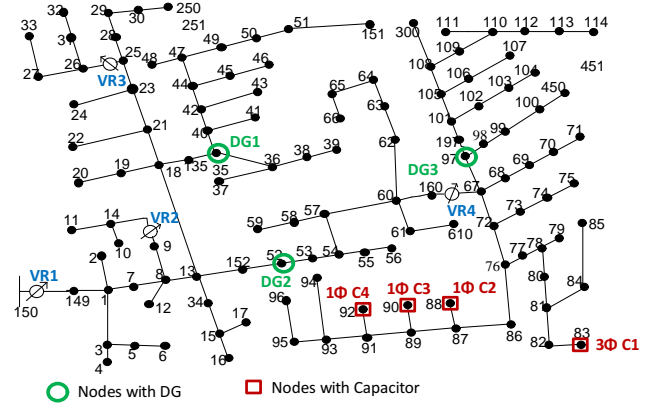}
\vspace{-0.4cm}
	\caption{IEEE 123-bus distribution test feeder.}
\vspace{-1.2cm}
	\label{fig:1}
\end{figure}

\vspace{-0.1cm}

\section{Results and Discussion}
The proposed VVO approach is validated using following three test feeders: IEEE 13-bus, IEEE 123-bus \cite{TestFeeder}, and PNNL 329-bus taxonomy feeder\cite{feeder329}. First, the proposed linear and nonlinear approximate power flow formulations are validated against the actual power flow solutions obtained using OpenDSS. Next, the proposed voltage-dependent load models are validated against equivalent ZIP load models. Finally, we demonstrate the proposed VVO approach using the aforementioned three test feeders. All simulations are done on MATLAB platform. Level-1 problem, modeled as MILP, is solved using CPLEX 12.7 and Level-2 problem, modeled as NLP, is solved using fmincon function in MATLAB optimization toolbox. A computer with core i7 3.41 GHz processor with 16 GB of RAM has been used for the simulations. The results obtained from MATLAB are validated against OpenDSS. 

IEEE-13 bus is a small highly loaded unbalanced distribution feeder operating at 4.16 kV making it a good candidate for testing VVO applications. This test feeder includes a three-phase and a single-phase capacitor bank and a voltage regulator at the substation. A PV with smart inverter of 575 kVA rated capacity is installed at node 671. IEEE-123 bus feeder also presents unbalanced loading conditions and several single-phase lines and loads with voltage drop problems making it a good candidate for demonstration of VVO application. There are four voltage regulators and four capacitor banks deployed along the feeder as shown in Fig. 1. The feeder is modified to include three DGs of capacity 345 kVA, 345 kVA, and 690 kVA at nodes 35, 52, and 97 respectively (see Fig. 1).  The 329-bus feeder is used to demonstrate the scalability of the proposed approach. Notice that 329-bus feeder includes 329 physical nodes and a total of 860 single-phase nodes. Compared to the state-of-art, this is a significantly large test system to demonstrate the coordinated control of all voltage control devices. The feeder includes one voltage regulator, one 600 kVAr three-phase capacitor bank, three 100 kVAr single-phase capacitor banks, and three DGs of capacity 23kVA, 57.5kVA and 115kVA (see Fig. 2).

\begin{figure*}[t]
\centering
\includegraphics[width=6.0 in ]{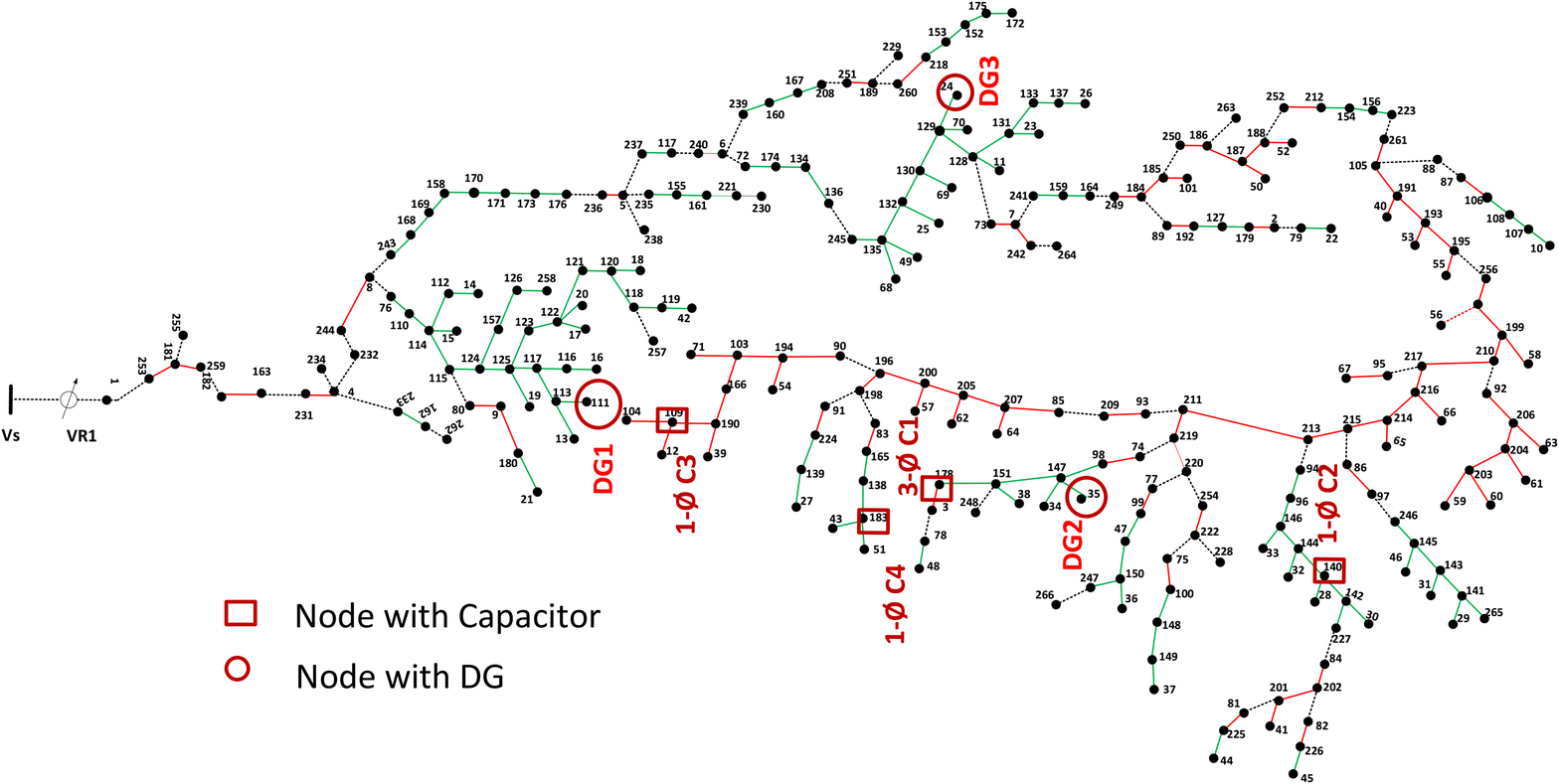}
\vspace{-0.6cm}
	\caption{PNNL 329-bus taxonomy distribution test feeder.}
\vspace{-0,5cm}
	\label{fig:1329}
\end{figure*}

Customer loads are assumed to have a CVR factor of 0.6 for active power and 3 for reactive power \cite{EPRI}. Note that the CVR values are arbitrary and can be easily adjusted based on the parameters for ZIP model of the load, if available, as detailed in Section III D.  To demonstrate the applicability of the proposed approach for different load mix, additional cases are simulated using a combination of residential and small and large commercial loads. The daily load and generation profiles are simulated in 15-min interval and are based on example profiles provided in OpenDSS (see Fig. 3). 

\begin{figure}[t]
\centering
\vspace{-0.2cm}
\includegraphics[width=3.3in, trim={0.3cm 0 0 0},clip]{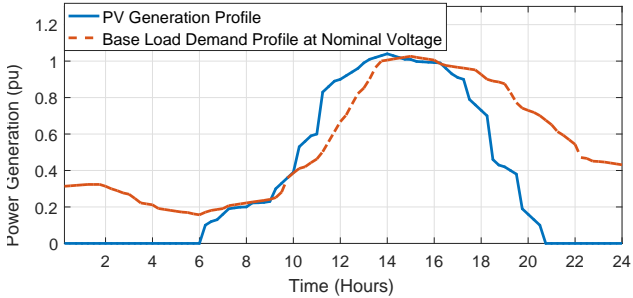}
\vspace{-0.4cm}
	\caption{Load demand and PV generation in 15-min interval.}
\vspace{-0.7cm}
	\label{fig:1lpv}
\end{figure}

\subsection{Verification of Approximate Power Flow Formulations}
\vspace{-0.2cm}
This section validates the proposed approximate power flow models. The results obtained from the proposed linear and nonlinear power flow models are compared with the power flow solution obtained using OpenDSS (see Table I). The largest errors in active and reactive power flow and bus voltages are reported for the three test feeders in Table I for different loading conditions. Note that the three-phase linear power model is sufficiently accurate in modeling power flow equations for an unbalanced system. Since, the losses are ignored in flow equations (equation (12)), the linear model incurs higher error in flow quantities ($P_{flow}$, and $Q_{flow}$). However, since the voltage drop due to flow quantities is included in linear model (equation (13)), the bus voltages are well approximated. Another key observation is an increase in error in $Q_{flow}$ vs. $P_{flow}$ for 123-bus and 329-bus feeders. The approximation errors in flow quantities using the linearized model depend upon relative values of line resistance and reactance. The line reactance is higher than the line resistance for these two feeders leading to more reactive power losses and hence a higher error in approximating $Q_{flow}$ quantities using linearized model. The nonlinear power flow model includes losses in the formulation and, therefore, results in lesser error in both flow quantities and bus voltages. The maximum error in bus voltages during peak load using linear and nonlinear models are: 0.0096 pu and 0.0025 pu for the 13-bus, 0.0074 pu and 0.0016 pu for the 123-bus, and 0.002 and 0.0002 pu for 329-bus systems, respectively. Note that 329-bus feeder is relatively more balanced and, therefore, incurs less error in voltages as compared to the rest of the two feeders. \\

	\begin{table}[t]
		\centering
		\caption{Comparison of Approximate Linear and Nonlinear Power Flow Formulations against OpenDSS Solutions}
\vspace{-0.2cm}        
		\label{singletable}
		\begin{tabular}{c|c|c|c|c}
			\toprule[0.4 mm]
			
			\multicolumn{5}{c}{Largest Error in Linear Power Flow wrt. OpenDSS Solutions }\\
			
			\toprule[0.4 mm]
			\hline
			{Test Feeder} & \% Loading & $P_{flow}(\%)$ &$Q_{flow}(\%)$& $V$(pu.)\\
			\hline
			\hline
			{IEEE 13 Bus}& 75\%& 5.1287& 4.938& 0.0075\\
            \hline
            {IEEE 13 Bus}& 100\%& 7.227& 6.442& 0.0096\\
            \hline
            {IEEE 123 Bus}& 75\%& 5.248& 9.502& 0.0054\\
            \hline
            {IEEE 123 Bus}& 100\%& 5.328& 11.313& 0.0074\\
            \hline
			{PNNL 329 Bus}& 75\%& 1.16& 6.9& 0.001\\
            \hline
            {PNNL 329 Bus}& 100\%& 1.55& 9.51& 0.002\\
			\toprule[0.4 mm]
			
			\multicolumn{5}{c}{Largest Error in Nonlinear Power Flow wrt. OpenDSS Solutions }\\
			\toprule[0.4 mm]
			{Test Feeder} & \% Loading & $P_{flow}(\%)$ &$Q_{flow}(\%)$& $V$(pu.)\\
			\hline
			\hline
			{IEEE 13 Bus}& 75\%& 0.2414& 1.668& 0.0015\\
            \hline
            {IEEE 13 Bus}& 100\%& 0.297& 2.034& 0.0025\\
            \hline
            {IEEE 123 Bus}& 75\%& 0.505& 2.58& 0.0014\\
            \hline
            {IEEE 123 Bus}& 100\%& 0.606& 3.88& 0.0016\\	
            \hline
            {PNNL 329 Bus}& 75\%& 0.3& 2.2& 0.0001\\
            \hline
            {PNNL 329 Bus}& 100\%& 0.6& 3.4& 0.0002\\		
            \toprule[0.4 mm]
		\end{tabular}
\vspace{-0.6cm}
	\end{table}
	
\vspace{-2.3cm}
\begin{table}[ht]
		\centering
		\caption{Maximum Error in approximating phase angle differences}
\vspace{-0.2cm}
		\label{singletableerr}
		\begin{tabular}{c|c|c|c}
			
			\toprule[0.4 mm]
			\hline
			{Test Feeder} & \% Load & error in $\delta_{ij}^{pq}$ & error in $\theta_i^{pq}$\\
			\hline
			\hline
			{IEEE 13 Bus}& 75\%& 1.8& 1.8\\
            \hline
            {IEEE 13 Bus}& 100\%& 2.1& 2.2\\
            \hline
            {IEEE 123 Bus}& 75\%&0.8 & 0.9\\
            \hline
            {IEEE 123 Bus}& 100\%&1.13 & 1.3\\
                       \hline
            {PNNL 329 Bus}& 75\%&0.5 & 0.55\\
            \hline
            {PNNL 329 Bus}& 100\%&0.9 & 1.05\\

			\toprule[0.4 mm]
			\end{tabular}
		\vspace{-0.5cm}
	\end{table}

The proposed nonlinear power flow formulation is based on two approximations: 1) difference between phase angles of node voltage ($\theta_i^{pq}$) is close to $120^\circ$, 2) difference between phase angles of branch currents ($\delta_{ij}^{pq}$) is close to those obtained using a constant impedance load model. These approximations are validated in Table II. The table reports largest deviation between actual quantities obtained using OpenDSS vs. the approximated ones used in this paper. As it can be seen, the largest error is less than $2^\circ$ for both voltage and current phase angle difference.

\begin{table}
		\centering
			\caption{Comparison of Approximate Nonlinear Power Flow Formulations against OpenDSS Solutions}
\vspace{-0.2cm}        
		\label{singletablec}
		\begin{tabular}{c|c|c|c|c}
			\toprule[0.4 mm]
			
			\multicolumn{5}{c}{Error in Nonlinear Power Flow wrt. OpenDSS for unbalanced case}\\
			
			\toprule[0.4 mm]
			\hline
			{Test Feeder} &  $V_{unbal} (\%)$ & $P_{flow}(\%)$ &$Q_{flow}(\%)$& $V$(pu.)\\
			\hline
			\hline
             {IEEE 123 Bus}& 3.2 & 0.61& 3.88& 0.002\\
            \hline
            {IEEE 123 Bus}& 5.7 & 0.68& 3.94 & 0.004 \\
            \hline
            {PNNL 329 Bus}& 2.8 & 0.52 & 3.74 & 0.0007 \\
            \hline
            {PNNL 329 Bus}& 4.5 & 0.95& 4.51& 0.0012\\
			\toprule[0.4 mm]
			\end{tabular}
\vspace{-1.5cm}
	\end{table}
The proposed nonlinear power flow is verified at heavily unbalanced loading for IEEE 123 node and 329 node system. The unbalanced in the system is created by increasing the load for one of the phase. The voltage unbalance  defined in (51) according to IEEE definition, is used to quantify the effect of load unbalance created in the system.  
\begin{equation} \label{eq51}
V_{unbalance} = \frac{max. deviation }{|V_{avg}|}*100
\end{equation}
The IEEE-123 node system has the inherent  apparent power unbalance of 23.2\%, which creates a maximum voltage unbalance of 3.2\%. Further, to produce more unbalance in the system the apparent power of phase A is increased which results in apparent power unbalance of 39.6\%. Due to increased power unbalance the maximum voltage unbalance in the system is 5.7\%. The maximum error in $P_{flow}$ ,$Q_{flow} $ and voltage is shown in Table III for IEEE-123 node system for both the test case. It is known that the 329 bus system is a balanced system. Hence, to generate unbalance in the system  the apparent power of phase B is increased which originates to apparent power  unbalanced of 40.72\% and 57.39\%. The effect of unbalanced power results in  voltage unbalanced  of 2.8\% and  4.5 \% respectively. The maximum error in  $P_{flow}$ ,$Q_{flow} $ and voltage with respect OpenDSS power flow results  is shown in Table III. It is required to mention that the nonlinear power flow is solved at flat start. According to ANSI C84.1 electric system can have the maximum voltage  unbalance of 3\%. Hence, Table III uphold the proposed power flow can be used for the heavily unbalanced system.
\vspace{-0.3cm}
\subsection{Validation of Proposed CVR-based Load Model}
 The proposed CVR-based voltage dependent load model derived in equations (21)-(22) is validated against equivalent ZIP load models detailed in equations (18)-(19). When accurately modeled, the CVR-based load model should require the same power demand as the equivalent ZIP load model for the acceptable range of operating voltages (0.95pu-1.05pu). Therefore, to validate the load models, the active and reactive power consumption for CVR-based load models are compared against the power consumption for ZIP load model for varying node voltages. ZIP models for residential, small commercial, and large commercial loads are used for validation. The ZIP coefficients for the different class of loads are obtained from \cite{Bokhari} and converted to CVR-based load model using equation (25) (also see Table IV).

The simulation case is detailed here. For each load class, the base active $p_{i,0}$ and  $q_{i,0}$ reactive power are taken as 100 kW and 100 kVAr, respectively. The voltage at the load point is varied from 0.95 to 1.05 pu. The active and reactive power demand for the two load models are shown in Fig. 4. It can be observed that for different load classes, the variation in power demand, both active and reactive, due to change in bus voltage are similar for both CVR-based load model and equivalent ZIP load model.
\begin{table}[t]

		\centering
		\caption{ZIP coefficients for different class of loads }
\vspace{-0.2cm}
		\label{singletablezip}
		\begin{tabular}{p{1.25cm}|p{0.4cm}|p{0.6cm}|p{0.4cm}|p{0.4cm}|p{0.73cm}|p{0.4cm}|p{0.5cm}|p{0.4cm}}
			\toprule[0.4 mm]
			\hline
			{Load Class} & $Z_p$ & $I_p$ & $P_p$ & $Z_q$ & $I_q$ &$Q_q$& $\hspace{-4pt}CVR_p$ & $\hspace{-4pt}CVR_q$\\
			\hline
			\hline
			{\hspace{-4pt}Residential}&  0.96 & -1.17 & 1.21 & 6.28 & -10.16 & 4.88&0.75&2.4\\
            \hline
            {\hspace{-4pt}Small Commercial}& 0.77 & -0.84 & 1.07 & 8.09 & -13.65 & 6.56&0.7&2.53\\
            \hline
            {\hspace{-4pt}Large Commercial}& 0.4 & -0.41 & 1.01 & 4.43 & -7.99 & 4.56&0.39&0.87\\
            			\toprule[0.4 mm]
			\end{tabular}
		\vspace{-0.1cm}
\end{table}

\begin{figure}[ht]
\centering
\vspace{-0.3cm}
\includegraphics[width=1.7in, trim={2cm 0 14cm 0},clip]{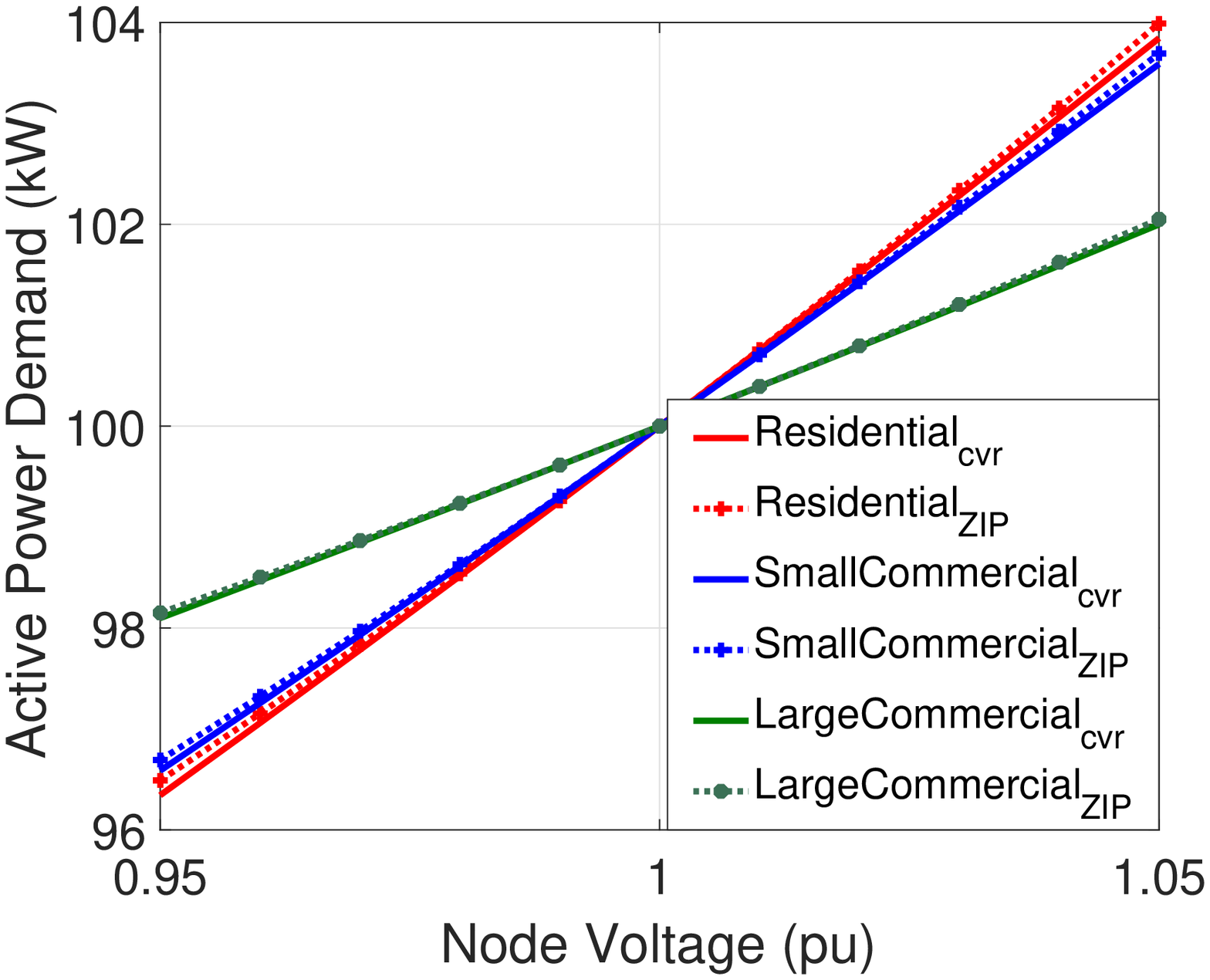}
\includegraphics[width=1.7in, trim={2cm 0 14cm 0},clip]{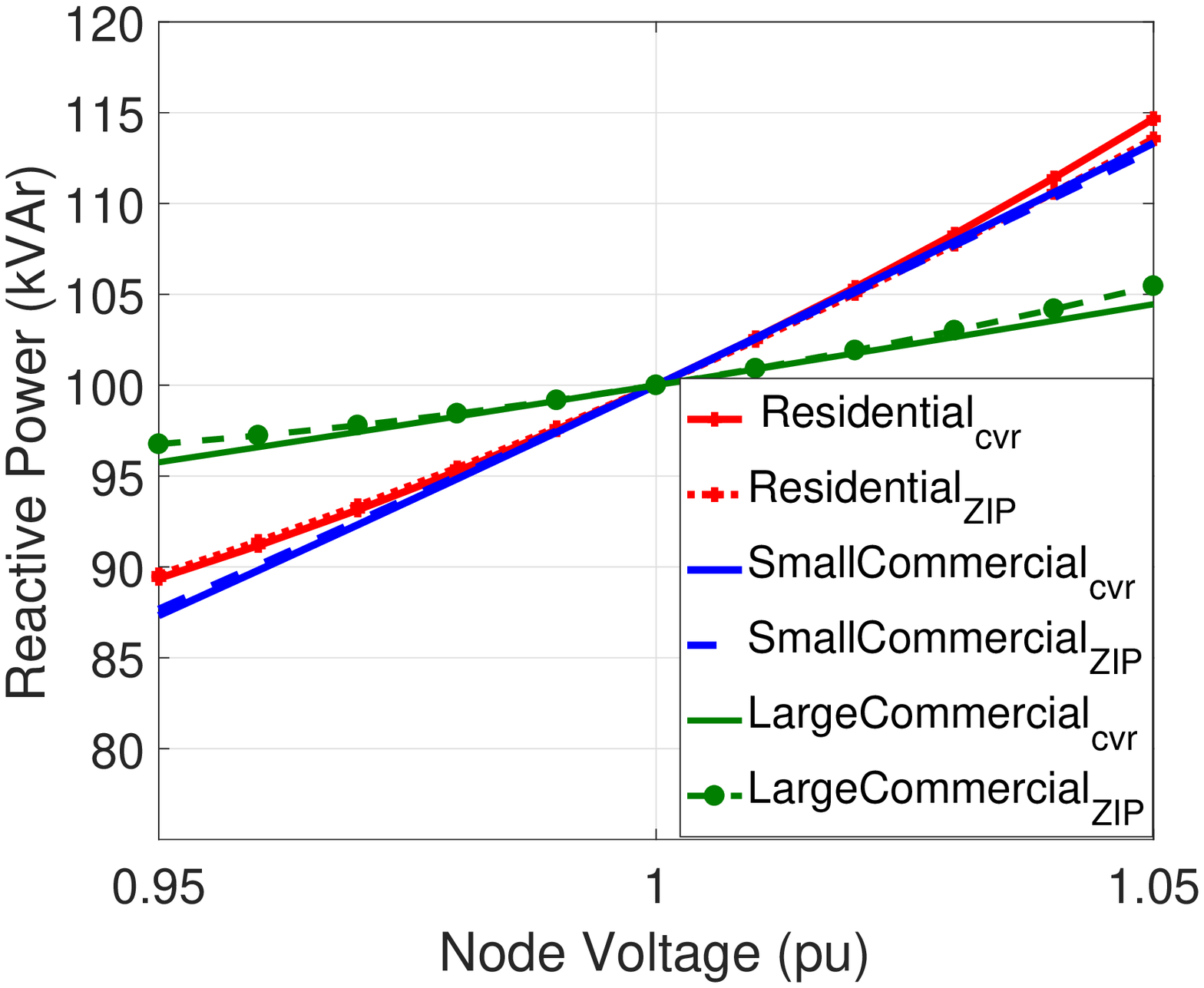}
\vspace{-0.3cm}
	\caption{\footnotesize Comparison of proposed load model with ZIP model: (a) Active power demand, (b) Reactive power demand.}
\vspace{-0.5cm}
	\label{fig:3}
\end{figure}

\subsection{CVR using proposed VVO approach}
The proposed bi-level VVO approach is validated using IEEE test feeders. The optimal control set points are obtained for both legacy and smart inverter control devices for the entire day. The results demonstrate that the proposed formulation ensures that feeder operates closer to minimum voltage range while not violating the voltage limit constraints and therefore, is effective in achieving CVR objectives.

\subsubsection{IEEE 13-bus test system}
The control variables for this feeder include a 32-step three-phase voltage regulator, a three-phase capacitor bank (Cap1), a single-phase capacitor bank on Phase C (Cap2), and one three-phase DG with smart inverter control. The model is simulated in 15-min interval for 1 day. The results obtained for the day during minimum and maximum loading conditions are shown in Table IV. As it can be seen that the feeder is unbalanced with a largest difference of around 0.24 MW during peak load condition. Based on the table, with the increase in load, regulator tap position changes from -13 at minimum load condition to 14 at maximum load condition. The three-phase capacitor bank is OFF for both load conditions, while single-phase capacitor in ON during peak load condition. For each phase, DG is supplying reactive power in order to maintain the feeder voltages within the ANSI limits except for phase C during maximum load condition. This is because the single-phase capacitor is ON and supplies the required reactive power. It should also be noted that the absorbed power supplied by DG increases with the loading. The average feeder voltage seen at both maximum and minimum loading conditions are close to 0.96 pu. 
The proposed VVO approach is, therefore, successful in maintaining feeder voltages close to minimum voltage limit, thus help extract the CVR benefits.

 The results obtained from the proposed approach are validated using OpenDSS. The optimal controls obtained from MATLAB for both maximum and minimum load conditions are implemented on 13-bus test feeder. The test feeder, with given statuses of voltage control devices, is solved using OpenDSS and substation power demand and minimum, maximum, and average node voltages are recorded (see Table V). It is observed that the system parameters obtained from MATLAB and OpenDSS closely match. This is expected given the accuracy of the proposed nonlinear power flow model.

\begin{table}[t]
		\centering
		\caption{Volt-VAR Optimization Results for IEEE 13-Bus Feeder ($CVR_p = 0.6$ and $CVR_q = 3$)}
\vspace{-0.3cm}
		\label{singletable13cvr}
		\begin{tabular}{c|c|c|c|c|c|c}
			\toprule[0.4 mm]
			IEEE-13 &\multicolumn{3}{|c|}{Minimum Load}& \multicolumn{3}{|c|}{Maximum Load}\\
            \hline
            \multicolumn{7}{c}{OPF solution from MATLAB}\\
            \hline
			Phase & A & B &C& A & B &C\\
			\hline
			{Regulator Tap}& -13 & -13 &-13& 14 & 14 &14\\
			\hline
			{Cap1 Status}& OFF & OFF & OFF & OFF & OFF &OFF\\
            \hline
            {Cap2 Status}& --- & --- &OFF& --- & --- &ON\\
            \hline
            {DG1 $q_{DG}^p$(MVAR)}& -0.04 & -0.13 &-0.03& -0.30 & -0.45 &0.12\\
            \hline
            \multicolumn{7}{c}{Optimal substation power flow and voltages using MATLAB}\\
            \hline
            {Load (MW)}& 0.143 & 0.132 &0.091& 0.866 & 0.622 &0.685\\
            \hline
            {Min. Voltage (pu)}& 0.955 & 0.956 &0.955& 0.95 & 0.955 &0.95\\
            \hline
            {Max. Voltage (pu)}& 0.97 & 0.97 &0.97& 1.03 & 1.03 &1.03\\
            \hline
            {Avg. Voltage (pu)}& 0.958 & 0.956 &0.958& 0.972 & 0.971 &0.972\\
			\hline
            \multicolumn{7}{c}{Validation of substation power flow and voltages using OpenDSS}\\
            \hline
			{Load (MW)}& 0.144 & 0.135  &0.095 & 0.87 & 0.625  &0.69\\
            \hline
            {Min. Voltage (pu)}& 0.955& 0.955& 0.955& 0.95& 0.953&0.95\\
            \hline
            {Max. Voltage (pu)}&0.97 &0.97 &0.97 & 1.03&1.03 &1.03\\
            \hline
            {Avg. Voltage (pu)}& 0.958&0.956 &0.958 &0.971&0.97 &0.971\\
			\toprule[0.4 mm]
\end{tabular}
		\vspace{-2cm}
	\end{table}

The optimal power consumption as recorded at the substation transformer after implementing the proposed VVO strategy for the day is shown in Fig. 5. The total power demand is compared with the case when VVO is not enabled. For this case, the capacitors and voltage regulators work in autonomous control mode while DG is operating at unity power factor. Except for peak demand duration, the proposed approach results in a reduction in net power demand. The largest reductions are seen at low loading condition.

 To further validate the proposed approach, we include additional test results with realistic load models for residential, commercial and large commercial loads. The ZIP coefficients details in Table IV are used to obtain CVR factors for each case with different load mix. The total load demand for minimum and maximum load conditions are reported in Table VI. It can be observed that the reduction in active power demand is lower for load mix with large commercial load as it shows less sensitivity to voltage.
\vspace{-0.7cm}
\begin{table}[ht]
\centering

		\caption{CVR for IEEE 13-Bus Feeder}
		\label{singletable2}
\begin{tabular}{p{2cm}|p{1cm}|p{1.3cm}|p{1cm}|p{1.3cm}}
\toprule[0.4 mm]
  {Load } & \multicolumn{2}{|c|}{Minimum Load} & \multicolumn{2}{|c}{Maximum Load} \\
  \cline{2-5}
  {Composition} &CVR& No CVR & CVR & No CVR \\
  \hline
  100\% R&  0.463 & 0.481& 2.096 &2.233 \\
  \hline
  70\% R, 30\% SC & 0.465 & 0.480  &2.098  & 2.234  \\
  \hline
  50\% R, 30\% SC, 20\% LC& 0.468 &0.475 &2.103 & 2.237 \\
  \hline
\end{tabular}
\footnotemark{R-Residential, SC-Small Commercial, LC-Large Commercial}
\vspace{-0.4cm}
\end{table}

\vspace{-0.7cm}
\begin{figure}[ht]
\centering
\includegraphics[width=3.6in]{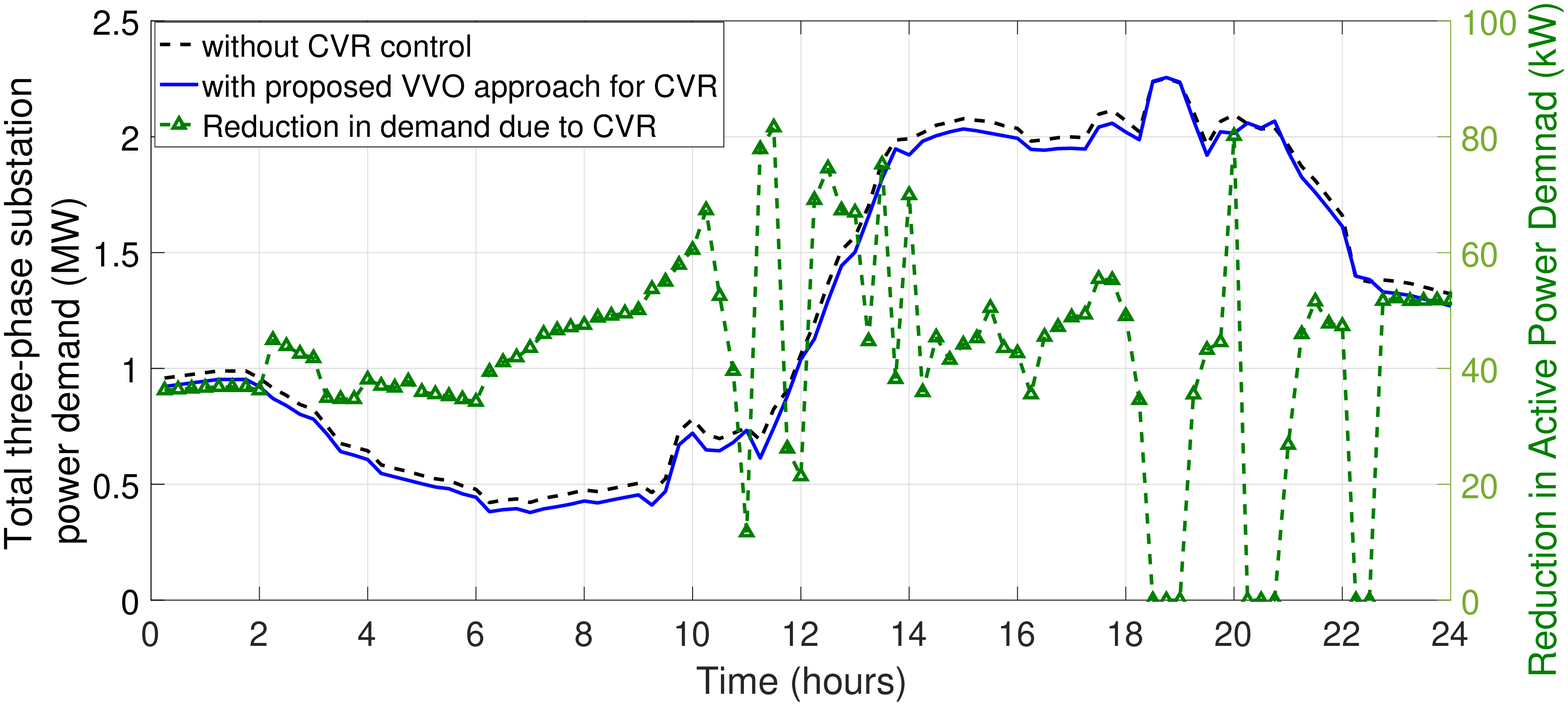}
\vspace{-0.8cm}
	\caption{IEEE-13 Bus CVR benefits observed using the proposed approach ($CVR_p = 0.6$ and $CVR_q = 3$).}
\vspace{-0.7cm}
	\label{fig:5}
\end{figure}

\vspace{1.8cm}

{\em Computational Complexity}
On an average on a dual core i7 3.41 GHz processor with 16 GB of RAM, the Level-1 solutions are obtained in less than 5 sec and Level-2 solutions are obtained within 9 sec for the 13-bus system. The largest time taken for solving Level-2 problem is 14 sec.

\begin{table}[t]

		\centering
		\caption{Volt-VAR Optimization Results for IEEE 123-Node Feeder ($CVR_p = 0.6$ and $CVR_q = 3$)}
\vspace{-0.1cm}
		\label{singletablecvr123}
		\begin{tabular}{c|c|c|c|c|c|c}
			\toprule[0.4 mm]		
			IEEE-123 &\multicolumn{3}{|c|}{Minimum Load}& \multicolumn{3}{|c|}{Maximum Load}\\
			\hline
			Phase & A & B &C& A & B &C\\
			\hline
            \multicolumn{7}{c}{OPF solution from MATLAB}\\
            \hline
			{Reg1 Tap}& -13 & -13 &-13& -8 & -8 &-8\\
			\hline
            {Reg2 Tap}& 0 & --- & ---& -2 &---& ---\\
			\hline
            {Reg3 Tap}& 1 & --- & 1& 7 & --- &2 \\
			\hline
            {Reg4 Tap}& 0 & 0 &0& 0 & 0 &0\\
			\hline
			{Cap1 Status}& OFF & OFF & OFF & ON & ON &ON\\
            \hline
            {Cap2 Status}& OFF & --- &---& ON & --- &---\\
            \hline
            {Cap3 Status}& --- & OFF & --- & --- & OFF &---\\
            \hline
            {Cap4 Status}& --- & --- & OFF & --- & ---&OFF\\
            \hline
            {DG1 $q_{DG}^p$(MVAR)}& -0.03 &0.045& 0.012 &-0.028 &0.03 &0.04\\
            \hline
            {DG2 $q_{DG}^p$(MVAR)}& 0.04 & -0.03 &0.03& -0.025& 0.039& -0.01\\
            \hline
            {DG3 $q_{DG}^p$(MVAR)}& -0.08 & -0.02 &-0.08& -0.09& 0.045& -0.09\\
            \hline
            \multicolumn{7}{c}{Optimal substation power flow and voltages using MATLAB}\\
            \hline
			{Load (MW)}& 0.20 & 0.13 &0.18& 0.99 & 0.78 &1.02\\
            \hline
            {Min. Voltage (pu)}& 0.955 & 0.955 &0.955& 0.951 & 0.953 &0.951\\
            \hline
            {Max. Voltage (pu)}& 0.965 & 0.965 &0.965& 0.995 & 0.995 &0.995\\
            \hline
            {Avg. Voltage (pu)}& 0.957 & 0.957 &0.958& 0.963 & 0.965 &0.966\\
            \hline
            \multicolumn{7}{c}{Validation of substation power flow and voltages using OpenDSS}\\
            \hline
			{Load (MW)}& 0.205 & 0.134 & 0.183 & 1.00 & 0.79 &1.024\\
            \hline
            {Min. Voltage (pu)}& 0.954 & 0.954 & 0.954 & 0.95&0.95 &0.95\\
            \hline
            {Max. Voltage (pu)}& 0.965&0.965 &0.965 &0.995 &0.995 &0.995\\
            \hline
            {Avg. Voltage (pu)}&0.956 &0.956 & 0.956& 0.96& 0.961&0.963\\
			\toprule[0.4 mm]
\end{tabular}
		\vspace{-1.0cm}
	\end{table}

\subsubsection{IEEE 123-bus test system}
Similarly, the proposed bi-level VVO approach is implemented using IEEE 123-node system for 1 day at 15-min interval. The results obtained from VVO for 123-node system are shown Table VII. The feeder is unbalanced with Phase B load being less than Phase A and C. The voltage regulator, Reg1, located at substation transformer (see Fig. 1), is at $-13$ tap for minimum load and $-8$ tap at maximum load. The voltage regulator, Reg4 is always at tap 0. Voltage regulators 2 and 3 are single and two-phase devices respectively and operate as optimization program instructs. Cap1 is a three-phase device and is OFF during minimum load and ON at maximum load condition and supplies required reactive power to maintain the voltage profile. Cap2, Cap3 and Cap4 are single phase devices and are ON/OFF depending upon the load demand. The DGs are located at three-phase nodes (see Fig. 1). Compared to minimum load condition, the reactive power demand or generation for DG1 and DG3 does not change significantly, except DG3. In contrast with minimum loading, DG3 is absorbing reactive power in Phase B during maximum load condition. Since Reg3 does not change the tap position, Phase B of DG3 adjusts the set points to account for the increase in load. Similarly, since there is no other VVC device between Reg1 and DG2, there is a drastic change in optimal DG behaviour between the two load conditions. The feeder voltage characteristics are also shown in Table VII. On an average the feeder operates close to minimum voltage limit, i.e. 0.96 pu, for both load conditions.

 The bi-level VVO approach is validated against OpenDSS. The optimal status of capacitor banks switch, voltage regulator tap, and reactive power reference to the DGs, obtained from MATLAB, are implemented on OpenDSS model for the 123-bus system. The substation power demand and feeder voltage characteristics obtained using MATLAB are validated against OpenDSS (see Table VIII). The system parameters obtained from MATLAB closely match to those obtained from OpenDSS thus validating the VVO model.

\begin{figure}[t]
\centering
\includegraphics[width=3.5in]{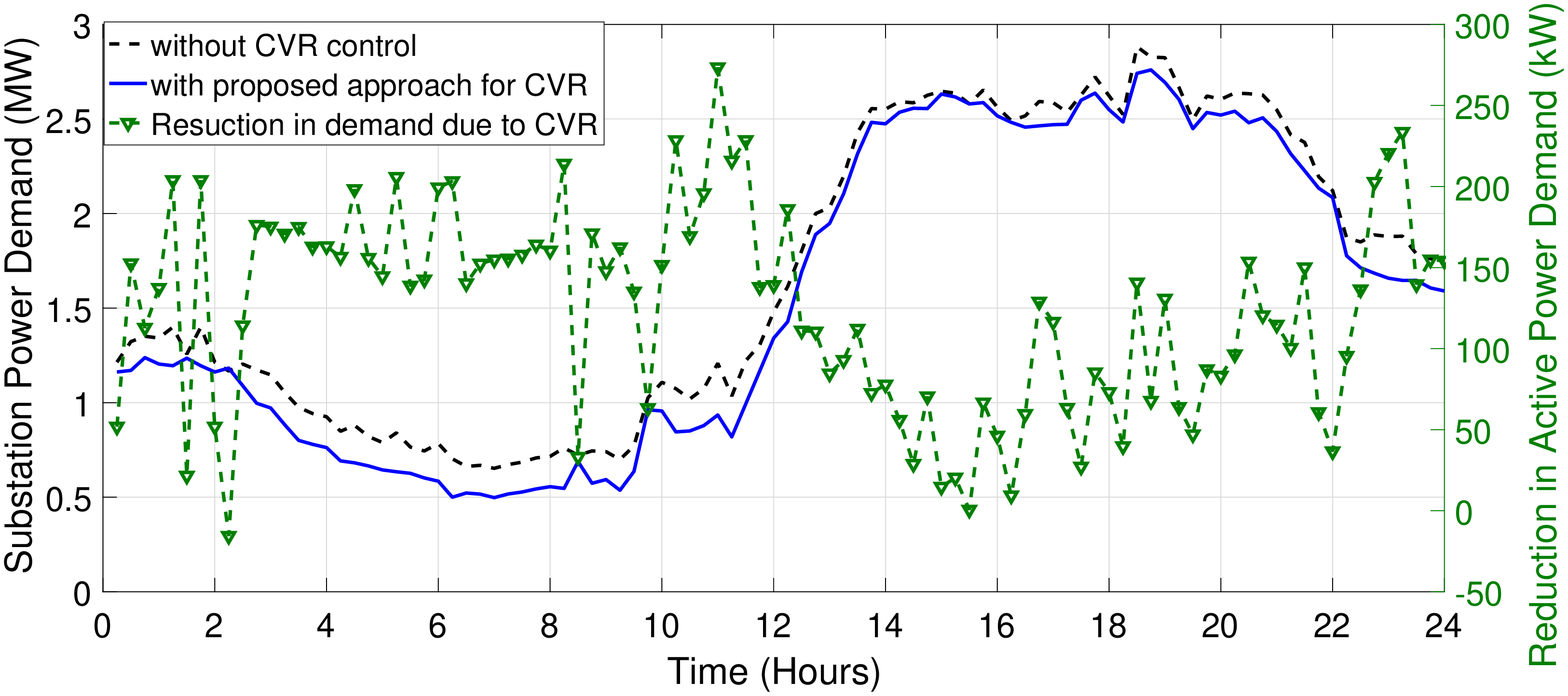}
\vspace{-0.7cm}
	\caption{\footnotesize IEEE-123 CVR Benefits Observed using the Proposed Approach ($CVR_p = 0.6$ and $CVR_q = 3$).}
\vspace{-1.4cm}
	\label{fig:7}
\end{figure}

\begin{table}[h]
\centering
\vspace{-0.25cm}

		\caption{CVR for IEEE 123-Bus Feeder}
		\label{singletable2mincvr}
\begin{tabular}{p{2cm}|p{1cm}|p{1.3cm}|p{1cm}|p{1.3cm}}
\toprule[0.4 mm]
    {Load } & \multicolumn{2}{|c|}{Minimum Load} & \multicolumn{2}{|c}{Maximum Load} \\
  \cline{2-5}
  {Composition} &CVR& No CVR & CVR & No CVR \\
  \hline
  100\% R&0.588 & 0.777 &2.726  &2.842 \\
  \hline
  70\% R, 30\% SC &0.588 & 0.776 &2.727  & 2.846  \\
  \hline
  50\% R, 30\% SC, 20\% LC&0.589  &0.748 &2.728 &  2.859\\
  \hline
\end{tabular}
\footnotemark{R-Residential, SC-Small Commercial, LC-Large Commercial}
\vspace{-0.2cm}
\end{table}

Finally, the CVR benefits obtained using the proposed approach are reported. The total three-phase substation load demand is compared to the case when VVO control is not enabled as shown in Fig. 6. On an average a reduction of around 150 kW is reported in net feeder active power demand. As expected the largest savings are reported during the minimum load condition.

 The proposed approach is further validated using ZIP load models for residential, commercial and large commercial loads. The ZIP coefficients detailed in Table IV are used to obtain CVR factors for each case with different load mix. The total feeder load demand for the minimum and maximum load condition are reported in Table VIII.

{\em Computational Complexity:}
On an average on a dual core i7 3.41 GHz processor with 16 GB of RAM, the Level-1 solutions are obtained in less than 5 sec for 123-bus system. For 123-bus, Level-1 solves 800 MILP equations with 1160 variables. On an average it takes 2-min to solve Level-2 problem for 123-bus system. The largest time taken for solving Level-2 problem for 123-bus system is 4 min. The Level-2 problem for 123-bus system solves 795 linear and 528 nonlinear equations (quadratic equalities) and with 1263 variables. The Level-1 and Level-2 solution times are within the 15-min control interval. Note that 123-bus test feeder represents a practical mid-size primary distribution circuit. The test feeder has 123 buses and a total of 267 single-phase nodes.

It should be noted that the Level-1 formulation scales well for larger feeders. This is because Level-1 solves an MILP that is relatively easier to solve even for a large set of constraints. The NLP problem in Level-2, however, is more difficult to scale for a large distribution system. In such cases, network reduction techniques are needed to represent the system with fewer equations \cite{NetRed}. In the following section, we demonstrate the scalability of the proposed approach using a 329-bus three-phase distribution feeder with the help of a simple network reduction technique.

%
\begin{table}[t]

		\centering
		\caption{Volt-VAR Optimization Results for 329-Node Feeder}
		\label{singletablecvr329}
		\begin{tabular}{c|c|c|c|c|c|c}
			\toprule[0.4 mm]	
			IEEE-329 &\multicolumn{3}{|c|}{Minimum Load}& \multicolumn{3}{|c}{Maximum Load}\\
			\hline
			Phase & A & B &C& A & B &C\\
            \hline
            \multicolumn{7}{c}{OPF solution from MATLAB}\\
			\hline
			{Reg1 Tap}& -6 & -6 &-6& 1 & 1 & 1\\
			\hline
           	{Cap1 Status}& OFF & OFF & OFF & ON & ON &ON\\
            \hline
            {Cap2 Status}& OFF & --- &---& OFF & --- &---\\
            \hline
            {Cap3 Status}& --- & OFF & --- & --- & OFF &---\\
            \hline
            {Cap4 Status}& --- & --- & OFF & --- & ---&OFF\\
            \hline
            {DG1 $q_{DG}^p$(MVAR)}& 0.02 &0.02& 0.02 &0.01 &0.01 &0.01\\
            \hline
            {DG2 $q_{DG}^p$(MVAR)}& -0.06 & -0.06 &-0.06& -0.03& -0.03& -0.03\\
            \hline
            {DG3 $q_{DG}^p$(MVAR)}& -0.11& -0.11 &-0.08& 0.055& 0.035& 0.022\\
            \hline
            \multicolumn{7}{c}{Optimal substation power flow and voltages using MATLAB}\\
            \hline
            {Load (MW)}& 0.444 & 0.459 &0.434& 2.86 & 2.97 &2.775\\
            \hline
            {Min. Voltage (pu)}& 0.958 & 0.958 &0.958 & 0.955 & 0.955 &0.955\\
            \hline
            {Max. Voltage (pu)}& 0.962 & 0.962 &0.962 & 1.0063 & 1.0063 & 1.0063\\
            \hline
            {Avg. Voltage (pu)}& 0.959 & 0.959 &0.959& 0.974 & 0.972 &0.976\\
            \hline
            \multicolumn{7}{c}{Validation of substation power flow and voltages using OpenDSS}\\
            \hline
            \hline
			{Load (MW)}& 0.445 &  0.462&0.438 & 2.87 & 2.98  &2.79\\
            \hline
            {Min. Voltage (pu)}& 0.958& 0.958& 0.958& 0.954& 0.953&0.954\\
            \hline
            {Max. Voltage (pu)}&0.962 & 0.962& 0.962&1.0063 &1.0063 &1.0063\\
            \hline
            {Avg. Voltage (pu)}&0.959 & 0.959& 0.959 & 0.971&0.97&0.973\\
			\toprule[0.4 mm]
\end{tabular}
		\vspace{-1.2cm}
	\end{table}

\subsubsection{329-bus PNNL Taxonomy Feeder}
The selected PNNL taxonomy feeder includes 329 buses, where, the number of nodes for phases A, B and C are 288, 298 and 274, respectively (total 860 single-phase nodes) (see Fig. 2). The proposed bi-level approach is implemented on 329-bus system. It is observed that Level-1 problem (MILP) takes on an average 20-sec. to solve, however, Level-2 problem (NLP) takes on an average 20-mins. Note that Level-2, for 329-bus system, solves for 4233 variables. In order to scale the Level-2 problem and to obtain a solution within 15-min interval, the 329-bus system is reduced using a simple network reduction technique. To reduce the network, we used the property of radial distribution feeders; the nodes that do not include branches, loads, or voltage control devices are combined using the equations for the series system for the corresponding branches. Using this method, the 329-bus system is reduced to a 184-bus system where, the number of nodes in phase A , B and C are 163, 171 and 156, respectively. After network reduction, the total number of variables for the Level-2 problem are reduced to 2415. Since network reduction is exact both models result in same power flow quantities. The Level-2 problem is implemented using 184-bus reduced network. The maximum computation time required to solve the Level-2 problem for the reduced network model is 9 mins.

\begin{figure}[t]
\centering
\includegraphics[width=3.4in, trim={1cm 0 0 0}, clip=true]{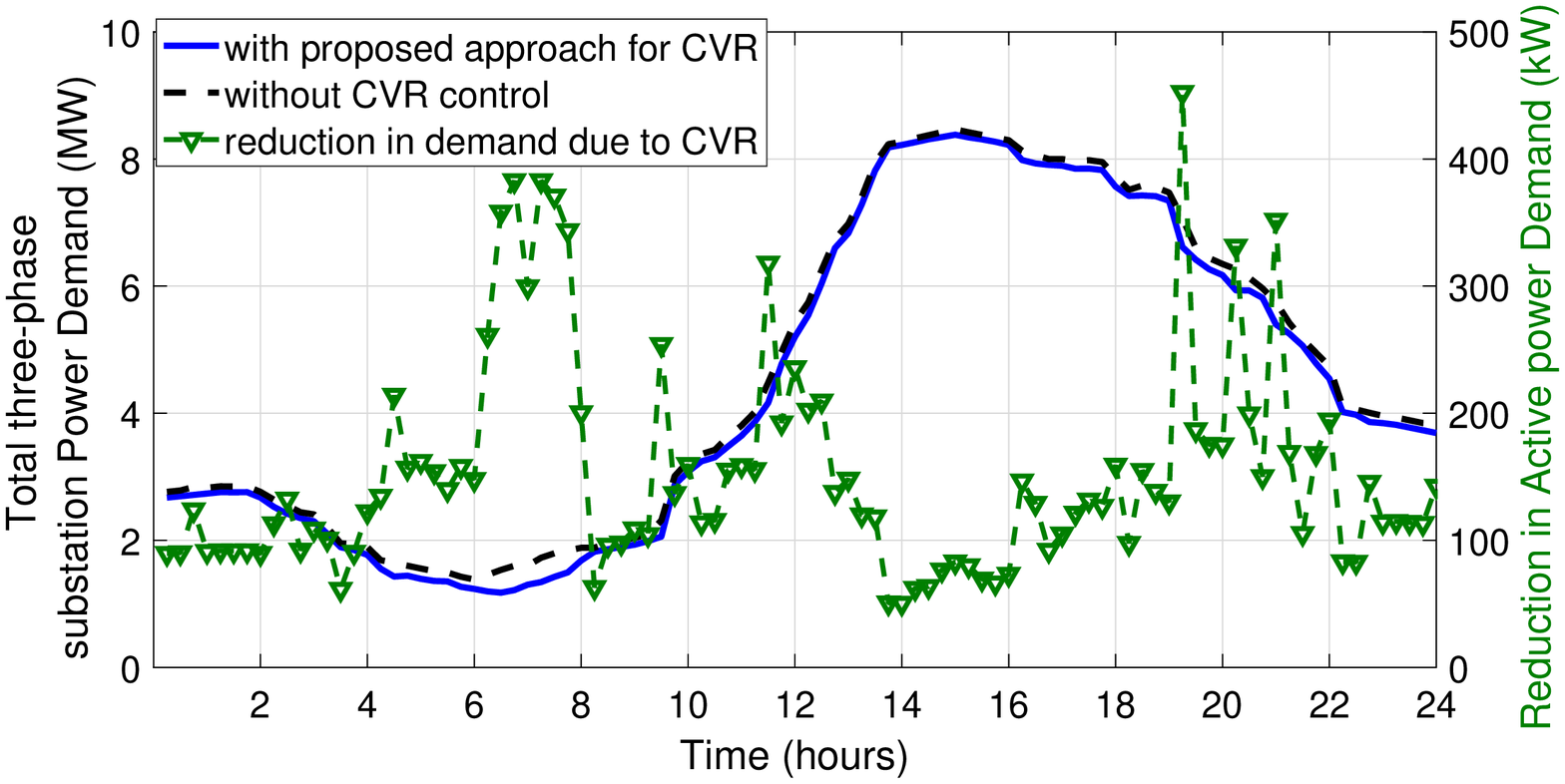}
\vspace{-0.4cm}
	\caption{\footnotesize 329-bus CVR Benefits Observed using the Proposed Approach ($CVR_p = 0.6$ and $CVR_q = 3$).}
\vspace{-0.5cm}
	\label{fig:9}
\end{figure}

The CVR results obtained for maximum and minimum load conditions are shown in Table IX. Note that the Level-1 problem is implemented using full 329-bus feeder and the Level-2 problem is implemented using reduced 184-bus feeder. As the load is closely balanced, the behavior of each phase is almost similar. The voltage regulator at the substation is at -6 tap for the minimum load and at 1 tap position for the maximum load condition. At minimum load, the three-phase as well as single-phase capacitor banks are OFF. However, at the maximum load condition, the three-phase capacitor is ON and single-phase capacitor  banks are OFF. The reactive power support from DG1 is same for all phases for both maximum and minimum load conditions. The minimum voltage for all the phases is at 0.958 pu at minimum load condition and at 0.955 pu at maximum load. The average voltage along the feeder is 0.959 and 0.972 at minimum and maximum load conditions, respectively. The substation power demand and feeder voltage characteristics obtained using MATLAB are validated against OpenDSS (see Table IX). The system parameters obtained from MATLAB closely match to those obtained from OpenDSS, validating the proposed VVO model.

The CVR benefits obtained using the proposed approach for 24-hour duration are reported in Fig. 7. The total three-phase substation load demand is compared to the case when VVO control is not enabled. On an average a reduction of around 200kW is reported in the net feeder active power demand.

\section{Conclusion}
This paper presents a VVO approach for CVR by coordinating the operation of distribution system's legacy voltage control devices and smart inverters. A bi-level VVO framework based on mathematical optimization techniques is proposed to efficiently handle the discrete and continuous control variables. The proposed approach solves OPF for a three-phase unbalanced electric power distribution system. The Level-1 solves a MILP problem to obtain control setpoints for both legacy devices and smart inverters using linear approximation for three-phase OPF. Next, Level-2 freezes the control for legacy devices and solves a NLP problem (with linear objective and quadratic constraints) to obtain a feasible and optimal solution by adjusting the setpoints for DG control using an approximate nonlinear OPF model. The approach is thoroughly validated using three test feeders, IEEE 13-bus, IEEE 123-bus, and PNNL 329-bus taxonomy feeders. The results demonstrate that: 1) the proposed power flow approximations are reasonably accurate, 2) the proposed approach successfully coordinates the operation of legacy and new devices for CVR benefits, and 3) both Level-1 and Level-2 solutions are computationally efficient for a realtime operation.

	
\bibliographystyle{ieeetr}
\vspace{-0.1cm}
\bibliography{references}

\end{document}